\documentclass[twoside,reqno,10pt]{amsart}
\usepackage{amsmath, amssymb} 
\usepackage{graphicx, color}
\usepackage{url}

\newcommand{\frdiffii}[3]{\ensuremath{\,_{#3}\mathbf{I}^{#1}_#2}}

\newcommand{\frdiffiix}[2]{\frdiffii{#1}{t}{#2}}

\newcommand{\frdiffir}[2]{\ensuremath{ \mathcal{D}^{#1}_{#2} }}
\newcommand{\wrighta}[3]{\ensuremath{ W \left(\left. #1,  #2 \right|  #3 \right) }}

\newcommand{\llim}[3]{\ensuremath{ \lim\limits_{ #1 \rightarrow #2} #3 }}
\newcommand{\fclass}[2]{\ensuremath{  \mathbb{#1}^{\, #2} }}
\newcommand{\erf}[1]{\ensuremath{ \mathrm{erf} \left( #1 \right)   }}

\newcommand{\epnt}{\; .}
\newcommand{\ecma}{\; ,}

\newtheorem{theorem}{Theorem}
\newtheorem{corollary}{Corollary}

\begin{document}

%     \begin{frontmatter}
     	
	\title{Anisotropic solutions of the time-fractional diffusion equation in multiple dimensions}

\author {Dimiter Prodanov }
\address{Correspondence: Environment, Health and Safety, IMEC vzw, Kapeldreef 75, 3001 Leuven, Belgium;  
	e-mail: Dimiter.Prodanov@imec.be, dimiterpp@gmail.com   }

%%%%%%%%%%%%%%
%  Abstract
%%%%%%%%%%%%%%
\begin{abstract}
	Anomalous diffusion phenomena are ubiquitous in complex media, such as biological tissues. 
	A wide class of sub-diffusive phenomena phenomena is described by the time-fractional diffusion equation. 
	The paper investigates the case of anisotropic fractional diffusion in the Euclidean space. 
	The solution of the fractional sub-diffusion equation can be expressed in terms of the Wright function and its spatial derivatives, parametrized by the directional unit vector (or alternatively a normal hyperplane).
	Moreover, the multidimensional case could be expressed as a transformation of the one-dimensional case. 

 
%   		
% Reviewers
% Jordan Hristov
%   		

 \end{abstract}
 %%%%%%%%%%%%%%%%%%%%%%%%%%
% 	\subsection*{Highlights}
% 	
% 	\begin{itemize}
% 		\item A concept of fractional-order velocity is defined. 
% 		\item Fundamental results about fractional-order velocity are established.
% 		\item It is demonstrated how derivative of a Holder function  can be regularized at the point where it diverges.
% 		\item Several cases of fractional Taylor expansions of Holder functions are investigated.
% 	\end{itemize}
  %%%%%%%%%%%%%%%%%%%%%%%%%%	
%  \end{frontmatter}
        	
     	\maketitle  

 %%%%%%%%%%%%%%%
 %  Sec
 %%%%%%%%%%%%%%
 \section{Introduction}\label{se:intro}
 
 Anomalous diffusion phenomena are ubiquitous in complex media, such as biological tissues \cite{Hoefling2013}, \cite{Metzler2017}. 
 The spatial complexity of the medium can impose geometrical constraints on transport processes on all length scales that can fundamentally alter the usual diffusion laws \cite{Metzler2004}. 
 
 A wide class of the observed sub-diffusive phenomena phenomena is described by the time-fractional diffusion equation. 
 This is so because the sub-diffusion can result from the continuous-time random walk (CTRW) model as an asymptotic limit \cite{Hilfer1995}, \cite{Compte1996}.
 In such case the order of the fractional derivative $0<\beta \leq 1$ describes the exponents of mean square displacement law: $<x^2>= \mathcal{O}\left(  t^\beta \right)  $ .
 Fractional diffusion models have been employed in hydrology describing well slow diffusion \cite{Meerschaert2006}; in advection-dispersion systems; or protein diffusion in the plasma membrane \cite{Kou2004}, porous biological tissues, lipid membranes  or the chromosomes (reviews in \cite{Pierro2018}, \cite{Sun2018}).

 Its one-dimensional version has been the subject of many papers  \cite{Kochubei1990},  \cite{Schneider1989},  \cite{Wyss1986}, \cite{Fujita1990}, \cite{Nigmatullin1986}. Notably, it is established that the fundamental solution on the real line is expressed by the  two-parameter  Wright function \cite{Gorenflo2000}.
 On the other hand, the multidimensional case has been addressed only assuming spherical symmetry of the solution, that is isotropy. 
 Schneider and Wyss obtained the solution of the time-fractional equation in terms of Fox functions \cite{Schneider1989}. They also showed that the Green’s function of fractional diffusion is a probability density.
 The fundamental solution of the  multidimensional spatially isotropic time-fractional equation was derived by Hanygad \cite{Hanygad2002} and later by Huang and  Liu \cite{Huang2005} using the conventional double Fourier-Laplace transform technique.  
 The resulting integral is, however, difficult to reverse-transform into the space-time domain.
 Recently, fractional multidimensional equations have been investigated by the combined use of Laplace, Fourier and Mellin transforms \cite{Ferreira2017}, \cite{Boyadjiev2017}.
 Ferreira and Viera also casted the problem in the language of complexified Clifford algebras.

 The present contribution instead uses only the Laplace transform techniques to derive the anisotropic fundamental solution. 
 The multidimensional solution can be also expressed by the Wright function and its spatial derivatives.
 The solution techniques is based on a the Dirac factorization procedure enabled by recasting the equation in the language of the Geometric algebra, that is the Clifford algebra over the Euclidean space \cite{Doran2003}.

 The manuscript is organized as follows: Section \ref{sec:formulation} introduces the transport problem, the main hypotheses and gives the solution on the real line. 
 Section \ref{sec:3dprob} starts by brief introduction of the Geometric Algebra of the Euclidean space and its extension to Geometric Calculus.  
 Section \ref{sec:green3d} derives the the anisotropic hyperplanar Green's function using the Dirac factorization technique. 
 In order to avoid detours into fractional calculus the main notation is given in Appendix \ref{sec:fract}. A brief overview of the properties of the Wright function is given in Appendix \ref{sec:wright}.
 
 %%%%%%%%%%%%
 %  Section
 %%%%%%%%%%%%
 \section{Formulation of the transport problem}\label{sec:formulation}
 
 In the Euclidean 3 dimensional space the transport problem can be formulated as
 \[
 \partial_{t}^{\beta+} c = D \nabla^2 c
 \] 
 where $c$ is the concentration of the species and  \textit{D} is a diffusion constant. 
 In addition, mass conservation law will be assumed to hold.
 Therefore, the concentration will be assumed to be normalized in the entire space
 \[
 \int_{ \fclass{R}{3}} c \ d x^3 = 1
 \]
 for all times.
 This constraint in turn implies that $ \llim{x}{\infty}{c (x)} = 0 $.
 
 This contribution focuses on a fractional derivative of the Caputo type because of its regular proprieties considering the terminals of integration. 
 Implicitly, it will be assumed that all of the considered unknown functions are of bounded variation either on the entire real line or on a suitable subinterval. 
 
 Furthermore,  functions from the kernels of the differ-integrals will also be excluded.    
 Notably, the kernels of Caputo and the Riemann-Liouville fractional derivatives are 
 \[
 Ker  [\frdiffir{\beta}{a+}]= \{f: \frdiffii{\beta}{t}{a+} f (t) = const \}
 \] 
 and 
 \[
 Ker  [ \partial_{t}^{\beta}  ]= \{ f (t) = const \}
 \]
 It should be noted that the kernels of the Caputo and  Riemann-Liouville fractional derivatives in general do not coincide since
 \[
 Ker  [\frdiffir{\beta}{a+}] \cap Ker  [ \partial_{t}^{\beta}  ]= \{0\}
 \]
 
 This contribution treats the $\nabla$ operator in a somehow different manner. 
 The transport problem in 3D will be solved using the methods of Geometric algebra.
 As this is not yet common knowledge in the physical sciences some introductory remarks are in order and are given in Sec. \ref{sec:geoalg}. 
 
 %%%%%%%%%%%%%
 %  Section
 %%%%%%%%%%%%%
 \subsection{The diffusion transport problem in one spatial dimension}\label{sec:1dprob}
 
 For simplicity we first consider the one-dimensional Cauchy problem :
 \[
 \partial_{t}^{\beta+} c  = D \, \partial_{xx} c 
 \] 
 with boundary conditions 
 \begin{flalign*}
 \llim{t}{0^-}{c (t,x)} &= 0, \quad \llim{x}{\infty}{c (t, x)} = 0 \\
 c(0, x) &= q (x), \quad t>0
 \end{flalign*}
 %%%%%%%%%%%%%%
 and initial condition $c(0,x) =\delta (x) $
 
 Since $\partial_{t}^{\beta}$ and $\partial_{x}$ commute the following factorization of the equation is in order 
 \[
 \left( \partial_{t}^{\beta/2} - \sqrt{D} \, \partial_{x}\right) \left( \partial_{t}^{\beta/2} + \sqrt{D} \,  \partial_{x}\right) c = 0
 \]
 as proposed by Oldham and Spanier \cite{Oldham1974}. 
 That is, giving due credit, this transformation can be called \textbf{Oldham -- Spanier factorization} in 1+1 dimension. 
 Therefore, the order of the equation can be reduced  and the solutions are given by the solutions of the system
 \begin{flalign}
 \partial_{t}^{\beta/2} c_{+} + \sqrt{D} \, \partial_{x} c_{+} = 0 \\ 
 \partial_{t}^{\beta/2} c_{-} - \sqrt{D} \, \partial_{x} c_{-} = 0
 \end{flalign}
 so that by linearity the general solution is
 \[
 c =  \lambda_{+} c_{+} + \lambda_{-} c_{-} ,
 \]
 where the sign index indicates the positive (resp. negative) spatial gradient.
 The general solution has been shown to be expressed by a special function of the Wright-type, called the M-Wright function.
 The calculation is straightforward but nevertheless it will be exhibited for completeness.

 The solution will be performed in the Laplace domain.
 The Laplace transform itself is introduced under the following notation
 \[
 \mathcal{L}_s : f(t) \div \hat{f} (s) = \int_{0}^{\infty} e^{ - s t} f(t) dt, \quad s \in \fclass{C}{}
 \]
 It is noteworthy that the Caputo derivative transforms as
 \[
 \mathcal{L}_s :	\partial_{t}^{\beta}  f(t) \mapsto s^\beta \hat{f} (s) - s^{\beta-1}  f (0^{+}) \, ,
 \]
 where the plus sign denotes passing to the limit as $t \rightarrow 0$.
 
 The full equation for the Green's function reads
 \[
 s^{\beta/2} \hat{G} -  s^{\beta/2 -1}  \delta(x)  + \sqrt{D} \,  \partial_{x} \hat{G} =  0
 \]
 Since $\delta(x)$ is a generalized function we will replace it with a smooth function
 $q_n(x)$, such that $\llim{n}{\infty}{q_n(x)} = \delta(x)$ point-wise, in the sense that 
 \[
 \llim{n}{\infty} {  \int_{-\infty}^{x} q_n(x) dx} = U(x), 
 \]
 the unit step function. 
 First, we look for the solution of the homogeneous equation 
 \[
 s^{\beta/2} \hat{G}_h  + \sqrt{D} \,  \partial_{x} \hat{G}_h =  0
 \]
 The solution is given straightforward by 
 $
 \hat{G}_h (s) = K e^{- x s^{\beta/2}/\sqrt{D} }
 $
 in the positive half-plane $x>0$ for the parameter \textit{K}, which possibly depends on 
 \textit{s}.
 A particular solution of the form $ \hat{G}_p (s) = p (x) e^{- x s^{\beta/2}/\sqrt{D} } $  will be  sought for further by variation of parameters.
 This results in the first-order ordinary differential equation (ODE)
 \[
 \partial_x p(x) \, \sqrt{D}\, {e}^{-\frac{x\, {s}^{\beta /2}}{\sqrt{D}}} -q_n(x)\, {s}^{\beta/2-1 }  =0
 \]
 This equation can be solved by direct integration as
 \[
 p(x)=  \frac{  {s}^{\beta/2-1 }  }{\sqrt{D}}\, \int_{-\infty}^{x} {\left. q_n( \xi)\, {  e}^{\frac{\xi\, {{s}^{\beta/2}}}{\sqrt{D}}}\, d \xi\right.}
 \]
 Passing to the limit for $\delta$ we obtain $p(x)= {s}^{\beta/2-1 } U(x)$.
 %%%%%%
 So  that
 \[
 \hat{G}_{+} (s) = \frac{e^{- x s^{\beta/2}/\sqrt{D} }}{\sqrt{D}\, s^{1-\beta/2}},  x\geq 0
 \]
 In a similar way by the reflection symmetry about the origin
 \[
 \hat{G}_{-} (s) = \frac{e^{ x s^{\beta/2}/\sqrt{D} }}{\sqrt{D}\, s^{1-\beta/2}},  x < 0
 \]
 so finally we take the union of the two solutions for the entire real line. 
 \[
 \hat{G} (s) = \hat{G}_{+} (s) \cup \hat{G}_{-} (s) =\frac{s^{ \beta/2 -1 }}{2} \frac{e^{ - |x| s^{\beta/2}/\sqrt{D} }}{\sqrt{D}}
 \]
 This solution is properly normalized for all admissible exponents $\beta$ since
 \[
 \mathcal{L}_s^{-1} \int_{-\infty}^{\infty} \hat{G} (s) dx = 1
 \]
 
 The inverse Laplace transform can be expressed as the Mainardi- Wright function (see Sec. \ref{sec:wright}):
 %%%%%%%%%%%%%%%%%%%%%%
 \begin{equation}\label{eq:green1d}
 G(x, t) = \frac{1}{2 \sqrt{D}\, t^{\beta/2}} M_{\beta/2} \left( \frac{|x|}{\sqrt{D} \, t^{\beta/2} }\right) = \frac{1}{2 \sqrt{D}\, t^{\beta/2}} \wrighta{-\frac{\beta}{2}}{1 -\frac{\beta}{2}}{-\frac{|x|}{\sqrt{D} \, t^{\beta/2} }}
 \end{equation}
 for $t>0$.
 %%%%%%%%%%%%%%%%

 %%%%%%%%%%%%%
 %  Section
 %%%%%%%%%%%%%
 \section{The transport problem in multiple dimensions}\label{sec:3dprob}

 For simplicity of the presentation let's consider first the problem in 3 spatial dimensions.
 
 %%%%%%%%%%%%
 %  Section
 %%%%%%%%%%%%
 \subsection{Primer on Geometric algebra}\label{sec:geoalg}
 
 The geometric algebra $\fclass{G}{3} \left( \fclass{C}{} \right)  $ is generated by the set of 3 orthonormal basis vectors $ E= \{e_{1}, e_{2}, e_{3} \}$
 for which the so-called \textit{geometric} product is defined with properties
 \begin{flalign}
 e_{1} e_{1}  &= e_{2} e_{2} = e_{3}  e_{3} = 1 \\
 e_{i} e_{j} &= -  e_{j} e_{i}, \quad i \neq j
 \end{flalign}
 %%%%%%%%%%
 An overview of the topic can be found, for example in the book of Doran and Lasenby \cite{Doran2003}.
 The geometric product of two vectors can be decomposed into a symmetrical scalar product and a  antisymmetrical wedge or exterior product
 \[
 a \, b = a \cdot b + a \wedge b \,  ,  \quad a \cdot b = b \cdot a \, , \quad  a \wedge b = -  b \wedge a
 \]
 Noteworthy, in \fclass{G}{3} there are zero divisors of the form
 \[
 \eta^{\pm} = \frac{1 \pm n}{2}, \quad n \cdot  n = 1 ,
 \]
 such that $ \eta^{+} \eta^{-} =0 $.
 Moreover, these elements are idempotents
 \[
 \eta^{\pm} \eta^{\pm}= \eta^{\pm}
 \]
 The bivector elements of the algebra $e_i e_j \equiv e_{ij}$ are isomorphic to the quaternion algebra and anti-commute. 
 The unique trivector element $e_{123} \equiv I$, that is the \textit{pseudoscalar}, commutes with all elements of the algebra and squares to $-1$ : $ I I =-1$.
 
 Geometric algebra has several advantages over the vector and vector analysis methods.
 On the first place, non-idempotent elements have inverses,  therefore, the division of vector elements is well-defined. 
 Secondly, the vector derivative can be treated as an element of the algebra.
 That is,  derivative operators can be multiplied safely and they also have inverses for suitable classes of functions.
 Lastly, the vector derivative is independent of the co-ordinates so that calculations can be performed in a coordinate-free manner.
 Altogether, this allows for much-greater flexibility than purely vector methods would allow.
 Moreover, by construction the vector methods are subset of geometric algebra methods.

 Consider the radius-vector $x=e_i x^j$ under the Einstein summation convention. 
 Then the vector derivative is defined as
 \[
 \nabla   := e^j \partial_{x_j} 
 \]
 where $ e^j =  \left( \partial_{x^j} x\right) ^{-1} $ are the elements of the dual basis, such that $e_i \cdot e^j =\delta_{ij} $.
 In \fclass{G}{3} they coincide with the elements of the usual basis \textit{E}. 
 For suitably smooth functions
 $
 \nabla \wedge \nabla =0
 $
 by commutativity of partial derivatives so that
 $
 \nabla \cdot \nabla = \nabla^2
 $.
 Then $\nabla r = p - q $, where the $p$ and $q$ are taken from the signature of the  algebra, so in the $n$-dimensional Euclidean case 
 $\nabla r = n $.

 %%%%%%%%%%%
 %  Section
 %%%%%%%%%%%
 \subsection{Dirac factorization technique for the diffusion equation}\label{sec:sol3d}
 
 Let's consider the transport problem in \textit{k} spatial dimensions in the Euclidean space. This corresponds to the Geometric algebra   $\fclass{G}{k} \left( \fclass{C}{} \right)  $.
 Since $\nabla$ is interpreted as a vector derivative then
 \begin{equation}\label{eq:diffga}
 \partial_{t}^{\beta} c = D \nabla \cdot \nabla c = D \nabla ^2 c
 \end{equation} 
 Therefore, the equation can be factorized as in the one dimensional case: 
 \[
 \left( \partial_{t}^{\beta/2} - \sqrt{D} \, \nabla \right) \left( \partial_{t}^{\beta/2} + \sqrt{D} \,  \nabla \right) c = 0
 \]
 We are looking for a homogeneous trial solution of the form
 $
 \hat{G} (s) = K e^ { r \cdot n \; s^{\beta/2} /\sqrt{D} }
 $
 which is parametrized by a direction vector $n =e_1 n_x + e_2 n_y + e_3 n_z$.
 Therefore,
 \[
 \left( s^{\beta/2} + \sqrt{D} \, \nabla \right) \left( s^{\beta/2} - \sqrt{D} \,  \nabla \right)  e^ { r \cdot n s^{\beta/2} /\sqrt{D} } = 0
 \]
 \[
 s^{\beta/2} \left( s^{\beta/2} + \sqrt{D} \, \nabla \right) \left( 1 - n  \right)   e^ { r \cdot n s^{\beta/2} /\sqrt{D} }  = 0
 \]
 so that finally
 \[
 s^{\beta} \left( 1 + n  \right)  \left( 1 - n  \right)   e^ { r \cdot n s^{\beta/2} /\sqrt{D} }  = 0
 \]
 Therefore,
 $
 1 - n^2=0
 $
 must hold and $n$ is a \textbf{constant unit vector}.
 Therefore, by duality $  n \wedge I (I^2)$ is the normal hyperplane  of constant phase.
 Hence, this provides an  diffusion analogue of the (hyper)plane wave solutions of the wave equation. 
 
 From this calculation it follows that in the Laplace domain the equation  
 \[
 \left( s^{\beta/2} - \sqrt{D} \,  \nabla \right) \hat{G}_n = s^{\beta/2} \left( 1 - n  \right) \hat{G}_n 
 \]
 holds  for the homogeneous fractional differential equation.
 Therefore, the grade of the equation can be lowered also in the time domain as follows
 \begin{equation}
 \left( \partial_{t}^{\beta/2} \pm \sqrt{D} \,  \nabla \right) G_n = 
 (1 \pm n) \, \partial_{t}^{\beta/2} G_n  = 2 \eta_{\pm} G_n
 \end{equation}
 In such way the relation to the idempotents and nilpotents of the algebra  becomes explicit.

 %%%%%%%%%%%
 %  Section
 %%%%%%%%%%%
 \section{The Hyperplanar Green's function}\label{sec:green3d}
 
 In multiple dimensions, the initial condition reads $c(0, r) =\delta_k (r) $, where $\delta_k$ indicates the \textit{k}-dimensional Delta function. 
 
 %%%%%%%%%%%%
 %% Section
 %%%%%%%%%%%%
 %\subsection{Explicit derivation of the hyperplanar Green's function }\label{sec:greenbox}
 
 In the Laplace domain in a similar way  the homogeneous equation reads
 \[
 s^{\beta} \hat{G}  - D \nabla^2 \hat{G} = 0
 \] 
 The homogeneous solution is given by
 $
 \hat{G}_h(s) = K   e^{ - n \cdot r  s^{\beta/2} /\sqrt{D} }
 $.
 Furthermore, the following normalization is in order for the positive-direction Green's function 
 \[
 \mathcal{L}^{-1}_t \int_{R_{+}^3} \hat{G}_h(s) \, d x^3 = \frac{1}{8} \Rightarrow   \int_{R_{+}^3} \hat{G}_h(s) \, d x^3  = \frac{1}{8 \, s}
 \]
 Therefore,
 \[
 \frac{1}{8 \, s}= \int_{R_{+}^3} \hat{G}(s) \, d x^3 = \left( \int_{0}^{\infty}  e^{ - n_i x_i  s^{\beta/2} /\sqrt{D}  }  d x_i \right)^3 = K \left( \frac{\sqrt{D}}{s^{\beta/2}}\right)^3
 \]
 by separability of the spatial variables, where we absorb some factors into \textit{K}. 
 Therefore, 
 \[
 \hat{G}(s) = \frac{s^{\frac{3 \beta}{2} -1} }{ (2 \sqrt{D})^3}  e^{ - n \cdot r  s^{\beta/2} /\sqrt{D} }
 \]
 The derivation technique can be extended further by induction. 
 Indeed, it can be noticed that the use of Cartesian coordinates ensure separability of the Laplace transformed solution in the spatial domain leading to 
 \[
 \frac{1}{2^k} =  \left( \int_{0}^{\infty}  e^{ - n_i x_i  s^{\beta/2} /\sqrt{D}  }  d x_i \right)^k
 \]
 Therefore, for $k$ dimensions the solution reads %it can be generalized
 %%%%%%%%%%%%%
 %  Eq
 %%%%%%%%%%%%%
 \begin{equation}
 \hat{G}_\mathbf{n}^k(s) = \frac{s^{\frac{k \beta}{2} -1} }{ (2 \sqrt{D})^k}  e^{ - \mathbf{n} \cdot r  s^{\beta/2} /\sqrt{D} }
 \end{equation}
 The time-domain solution is given then by Laplace transform inversion as the Wright function
 \[
 G_\mathbf{n}^k (r, t) = \frac{1}{(2 \sqrt{D})^k  t^{k \beta/2} }  
 \wrighta{-\frac{ \beta}{2}}{1 -  k \frac{ \beta}{2}}{-\frac{|\mathbf{n} \cdot r|}{ t^{\beta /2} \sqrt{D}}}
 \]
 which can be continued for all  directions by reflections across the origin .

 It is noteworthy that from  Eq \ref{eq:wrightdiff} it immediately follows that 
 %%%%%%%%%%%%%%
 %  Prop
 %%%%%%%%%%%%%% 
 \begin{theorem}[Hyperplanar Green function]\label{prop:raygreen1}
 	The anisotropic Green's function in the k-dimensional Euclidean space is given by 
 	\[
 	G_\mathbf{n}^{k}(r, t) = \left(- \frac{1}{2  }\right)^k   \; \left( \mathbf{n} \cdot \nabla \right) ^k  \wrighta{  -\frac{\beta}{2}}{ 1}{ - \frac{|\mathbf{n} \cdot r|}{ t^{\beta /2} \sqrt{D} }},
 	\]
 	where $\mathbf{n}$ is a constant unit vector specifying the direction.
 \end{theorem}
 %%%%%%%%%%%%%%
 \begin{corollary}
 	In terms of the M-Wright function the solution is given by
 \[
 G_\mathbf{n}^{k}(r, t) = \left(- \frac{1}{2  }\right)^{k} \frac{1}{\sqrt{D}  t^{\beta/2} }   \; \left( \mathbf{n} \cdot \nabla \right) ^{(k-1)}  \wrighta{  -\frac{\beta}{2}}{ 1 -\frac{\beta}{2}}{ - \frac{|\mathbf{n} \cdot r|}{ t^{\beta /2} \sqrt{D} }}
 \]
 \end{corollary}
 %%%%%%%%%%%%%%%%
 
 The function 
 \[
 W_I (z, \beta):=\wrighta{  -\frac{\beta}{2}}{ 1}{ - z}
 \]
 can be considered as a scale-invariant primitive of the k-dimensional Green's function.
 The values of the function can be calculated by the following Theorem:

 %%%%%%%%%%%%
 %  Theorem
 %%%%%%%%%%%%
 \begin{theorem}\label{th:imwright}
 	 The Integral Wright function can be represented as
 	 \[
 	 W_I(a, z) =\frac{1}{2} - \frac{1}{\pi} \int_{0}^{\infty} K( a,-z, r)  dr, \quad a <0, z \geq 0
 	 \]
 	  \[
 	 W_I(a, z) =\frac{1}{2} + \frac{1}{\pi} \int_{0}^{\infty} K( a,z, r)  dr, \quad a <0, z<0
 	 \]
 	 \[
 	K(a, z, r) = \frac{{{  e}^{\frac{-\cos{\left( \ensuremath{\pi}  a\right) } z}{{{r}^{a}}}-r}}\, \sin{\left( \frac{\sin{\left( \ensuremath{\pi}  a\right) } z}{{{r}^{a}}}\right) }}{r}
 	\]
 
 \end{theorem}
 %%%%%%%%%%%%%%%%

\begin{figure}
	\centering
	\includegraphics[width=0.8\linewidth]{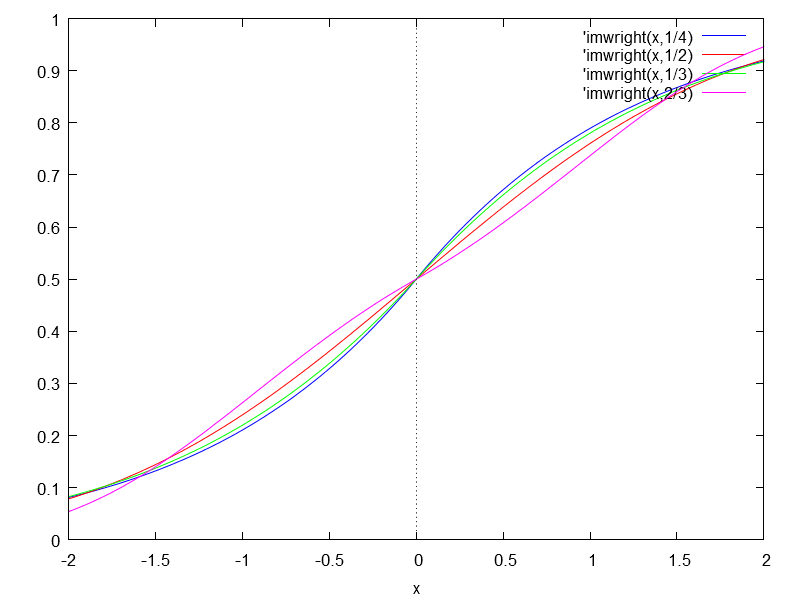}
	\caption{Plots of the Integral Wright function}
	The function's values are plotted for $\beta =\{1/4, 1/2, 1/3, 2/3 \} $ in blue, red, green and magenta, respectively. 
	\label{fig:imwright}
\end{figure}

Its plot is presented in Fig. \ref{fig:imwright}.
%%%%%%%%%%%%
%  Fig. 1
%%%%%%%%%%%%%

Finally, considering the usual applications,
for 3 spatial dimensions the solution can be specialized to
\[
G_\mathbf{n}^{3}(r, t) = \frac{1}{(2 \sqrt{D})^3  t^{3 \beta/2} }  
\wrighta{-\frac{ \beta}{2}}{1 -  3 \frac{ \beta}{2}}{-\frac{|\mathbf{n} \cdot r|}{ t^{\beta /2} \sqrt{D}}}
\]

 \section{Discussion}\label{sec:disc}
 
 The present work demonstrated that solutions of the $n$-dimensional fractional diffusion equation can be built from the solutions of the one-dimensional equation on the real line  by directional differentiation
 The solution method avoids the necessity of a second integral transform as originally performed by Mainardi et al. \cite{Mainardi2001}.
 Therefore, it has less assumptions about the spatial behavior of the Green's function or the spatial symmetries of the solution.   
 Notably, the existence of the solution of the entire real line need not be proven.
 This may be of advantage if one seeks solution in compartmentalized spatial domains. 
 
 Furthermore, the anisotropic hyperplanar solutions provide an interesting analogy with the plane wave solutions to the wave-equation, and could be extended in a generalized diffusion-wave formalism. 
 Such task, however goes beyond the scope of the present work and will be pursued elsewhere. 
 
 %%%%%%%%%%%
 %  Acknowledgments
 %%%%%%%%%
 \section*{Acknowledgments}
 The work has been supported in part by a grant from Research Fund - Flanders (FWO), contract number  VS.097.16N. %Graphs are prepared with the computer algebra system Maxima.
 
 %%%%%%%%%%%%
 %  Appendix
 %%%%%%%%%%%%
 \appendix
 
  \section{Appendices}\label{sec:appendice}
 %%%%%%%%%%%
 % Section
 %%%%%%%%%%%
 \subsection{Fractional Differ-integrals}\label{sec:fract}
 
 The left Riemann-Liouville differ-integral of order $\beta \geq 0$ (Samko et al. \cite{Samko1993} [p. 33]) is defined as
 \[
 \frdiffiix{\beta}{ a+ }  f (t) := 
 \frac{1}{\Gamma(\beta)} \int_{a}^{t}   f \left( \xi \right)  \left( t- \xi \right)^{\beta -1}d \xi 
 \]
 while the right integral is defined as
 \[
 \frdiffiix{\beta}{ -a }  f (x) = 
 \frac{1}{\Gamma(\beta)} \int_{t}^{a}   f \left( \xi \right)  \left( \xi-t \right)^{\beta -1}d \xi 
 \]
 where $\Gamma(x) $ is the Euler's Gamma function. 
 Here we depart from the definition of the right integral and complexify it in a naive manner as:
 \[
 \frdiffiix{\beta}{ -a }  f (x) :=  
 \frac{(-1)^{\beta-1}}{\Gamma(\beta)} \int_{t}^{a}   f \left( \xi \right)  \left( \xi-t \right)^{\beta -1}d \xi 
 \]
 which amounts to keeping the   order of variables as in the left integral.
 The left (resp. right) Riemann-Liouville (R-L) fractional derivatives are defined as the expressions (Samko et al. \cite{Samko1993}[p. 35]):
 \begin{flalign}
 \mathcal{D}_{a+}^{\beta} f (t)  & := \frac{d}{dt} \frdiffiix{1-\beta}{ a+ }  f (t)  = \frac{1}{\Gamma(1- \beta)}  \frac{d}{dt}  \int_{a}^{t}\frac{  f ( \xi ) }{{\left(  t - \xi \right) }^{\beta }}d \xi  \\
 \mathcal{D}_{-a}^{\beta} f (x) & := \frac{d}{dx}  \frdiffiix{1-\beta}{ -a }  f (x) = \frac{(-1)^{\beta}}{\Gamma(1- \beta)} \frac{d}{dt}  \int_{t}^{a}\frac{  f ( \xi )  }{{\left( \xi - t \right) }^{\beta }}d \xi
 \end{flalign}
 %%%%%%%%%%%%%%%%
 
 Under the same naming conventions, the fractional derivative in Caputo's sense  are defined as 
 %%%%%%%%%%%%%%%
 %  Definition
 %%%%%%%%%%%%%%%
 \begin{flalign}
 \partial_{t}^{\beta+}  f(t) & :=  \frdiffir{\beta}{a+} \left[  f-f(a) \right] (t) 
 = \frac{1}{\Gamma(1- \beta)}   \int_{a}^{t}\frac{  f^\prime ( \xi )   }{{\left( t- \xi \right) }^{\beta }}d \xi    \ecma \\
 \partial_{t}^{\beta-} f(t) & := (-1)^{\beta}  \frdiffir{\beta}{-a} \left[  f - f(a)   \right] (t)   =  \frac{(-1)^{\beta}}{\Gamma(1- \beta)}    \int_{t}^{a}\frac{  f^\prime ( \xi )   }{{\left( \xi - t \right) }^{\beta }}d \xi    \epnt
 \end{flalign}
 %%%%%%%%%%%%%%%%%%
 Wherever suitable the coordinate indices are skipped from notation.

 The left (resp. right) R-L derivative of a function $f$ exists for functions  representable by fractional integrals of order $\alpha$ of some Lebesgue-integrable function.
 This is the spirit of the definition of Samko et al. (\cite{Samko1993}, Definition 2.3, p. 43) for the appropriate functional spaces: 
 \begin{flalign*}
 \mathcal{I}^{\alpha}_{a, +} (L^1)& := \left\lbrace f: 
 \frdiffiix{ \alpha}{ a+ }  f (x) \in AC([a, b]), f  \in L^1 ([a,b]) , x \in [a, b] \right\rbrace \ecma  \\
 \mathcal{I}^{\alpha}_{a, -} (L^1)& := \left\lbrace f: 
 \frdiffiix{ \alpha}{ -a }  f (x) \in AC([a, b]), f  \in L^1 ([a,b]) , x \in [a, b] \right\rbrace 
 \end{flalign*}     
 Here $AC$ denotes absolute continuity on an interval in the conventional sense. 
 
 %%%%%%%%%%%
 %  Section
 %%%%%%%%%%%
 \subsection{Special Functions}\label{sec:specfun}
 
 %%%%%%%%%%%
 %  Section
 %%%%%%%%%%%
 \subsubsection{Bessel and Wright Functions}\label{sec:bessel}

 The Bessel $J$ function is represented by the integral
 \[
 J_\nu (z) = \frac{z^\nu}{(2 \pi)^{\nu+1}} \omega_{2 \nu} \int_{0}^{\pi} e^{- i z \cos{\theta}} \sin^{2\nu}{\theta} d \theta, \quad \omega_n = \frac{2 \ \pi^{n+\frac{1}{2}} }{\Gamma \left( n+ \frac{1}{2}\right) }
 \]
 For half-integer values the Bessel J function can be expressed by elementary functions. 
 For example
 \[
 J_{1/2} (z) = \sqrt{ \frac{ 2}{\pi}} \frac{ \sin{z}}{\sqrt{z}}
 \]

 %%%%%%%%%%%
 %  Section
 %%%%%%%%%%%
 \subsubsection{The Wright  function}\label{sec:wright}
 The  function \wrighta{\lambda}{\mu}{z},   named after E. M. Wright, is defined as the infinite series
 \begin{equation}
 \label{eq:wright}
 \wrighta{\lambda}{\mu}{z}:= %\Gamma(\mu)\HGFa{-}{-}{-}{  (\mu, \lambda)}{ z} = 
 \sum_{k=0}^\infty \frac{z^k }{k! \, \Gamma(\lambda k+\mu) },\quad \lambda >-1,\ \mu \in \fclass{C}{},
 \end{equation}
 \wrighta{\lambda}{\mu}{z} is an entire function of z.
 The summation is carried out with steps, such that $\lambda k+\mu \neq 0$.
 The function is related to the Bessel functions $J_\nu(z)$ and $I_\nu(z)$ as
 \begin{flalign*}\label{eq:besselj}
 \wrighta {1}{\nu+1} {- \frac{1} {4} z^2 } =\left( \frac{z}{2}\right)^{-\nu} J_\nu(z) \\
 \wrighta {1}{\nu+1} {  \frac{1} {4} z^2 } =\left( \frac{z}{2}\right)^{-\nu} I_\nu(z)
 \end{flalign*}
 and is sometimes called generalized Bessel function.
 A recent survey about the properties of the function can be found in \cite{Mainardi2010}.
 
\begin{figure}
 	\centering
 	\includegraphics[width=0.5\linewidth]{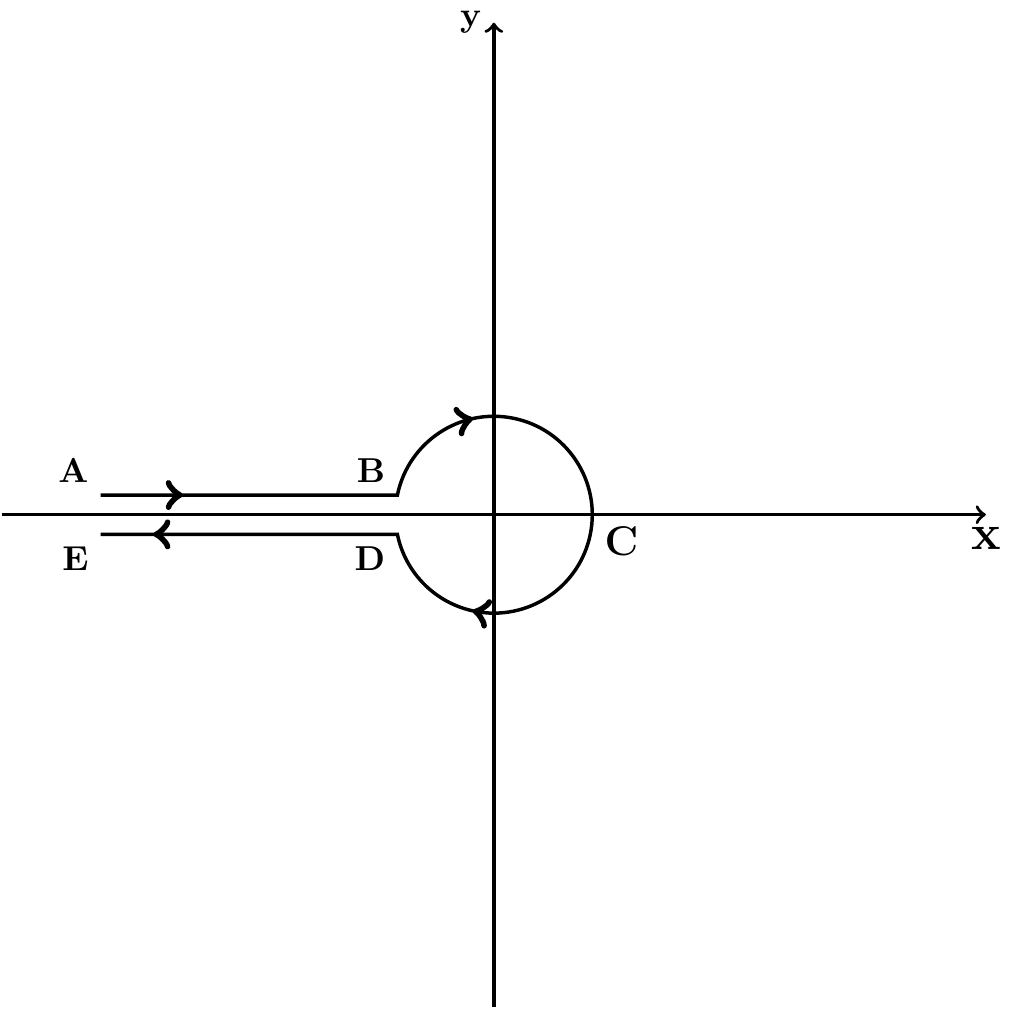}
 	\caption{Hankel contour}
 	\label{fig:hankel}
 \end{figure}

 The integral representation of the Wright function is noteworthy because it can be used for numerical calculations
 \begin{equation}\label{eq:intrep}
 \wrighta{\lambda}{\mu}{z}=\frac{1} {2\pi i}\int_{{\mathrm Ha}}
 e^{\zeta+z \zeta^{-\lambda}} \zeta^{-\mu}\, d\zeta, \quad \lambda >-1,\ \mu \in \fclass{C}{}
 \end{equation}
 where ${\mathrm Ha}$ denotes the Hankel contour in the complex $\zeta$-plane with a cut along the negative real semi-axis $\arg \zeta =\pi$. 
 	The contour is depicted in Fig. \ref{fig:hankel}.
 
 %%%%%%%%%%%
 %  Figure
 %%%%%%%%%%%%
 %%%%%%%%%%%%%%
 Furthermore,
 \begin{equation}\label{eq:wrightdiff}
 \frac{d}{dz} \wrighta{\lambda}{\mu}{z} =  \wrighta{\lambda}{\lambda+ \mu}{z}
 \end{equation}
 and formally
 \begin{equation}\label{eq:wrightint}
 \int \wrighta{\lambda}{\mu}{z} dz = \wrighta{\lambda}{\mu- \lambda}{z}   + C
 \end{equation}
 
 Mainardi introduces a specialization of the Wright function, which is called here the M-Wright function, which is important in the
 applications to fractional transport problems \cite{Mainardi2001}.
 \[
 M_{\nu}( z) := W(- \nu,1- \nu| -z)
 \]
 
 The M- Wright function can be calculated from its integral representation
 \begin{equation}\label{eq:intrepm}
 M_{\nu}( z) =\frac{1} {2\pi i}\int_{{\mathrm Ha}}
 \frac{e^{\zeta - z \zeta^{\nu}}} {\zeta^{1- \nu}}\, d\zeta 
 \end{equation}
 In particular, the following formula can be used for real arguments \cite{Luchko2008}
 %%%%%%%%%%%%%%%%%%%%
 \begin{flalign}
 M_{\nu}( z) &= \frac{1}{\pi \nu} \int_{0}^{\infty} K_\nu (r, z)\, dr \\
 K_\nu(r, z) &= {{  e}^{  -u z \cos{\pi\nu } -{{u}^{{1 / a}}}}}\, \sin{   \left( u z \, \sin{\pi\nu}- \pi \nu \right) }
 \end{flalign}
 %%%%%%%%%%%%%%%%
 The advantage here is that the integral kernel is not singular and allows for efficient computation. 
 Special cases of the M-Wright function are
 \[
 \begin{array}{c|c}
 %\hline 
 \nu & M_\nu (z) \\ 
 \hline 
 +	0 &  e^{-z} \\
 1/2 & \frac{1}{\sqrt{\pi}} e^{-z^2 /4} \\
 1/3 &  \sqrt[3]{3^2} Ai \left( z/ \sqrt[3]{3}\right)  \\
 \hline 
 \end{array} 
 \]

 %%%%%%%%%%%
 %  Section
 %%%%%%%%%%%
% \subsection{Laplace transforms of the Wright function}\label{sec:laplace}
 The Laplace transform image of a function $f$ will be denoted with $\hat{f}$ and the Laplace variable will be denoted by $s$ as
 \[
 \mathcal{L}_s : f(t) \mapsto \hat{f} (s)
 \]
 
 The M-Wright and Wright functions have some useful Laplace transform pairs \cite{Mikusinski1959}, \cite{Stankovic1970}, \cite{Gorenflo1999}
 %%%%%%%%%%%%%%%
 \begin{flalign}\label{eq:laptrsnfpairm}
 \frac{\nu k  }{t^{\nu+1}}  M_{\nu}\left(  k t^{-\nu}\right) \div  \exp(-k s^\nu),\ 0<\nu<1,\ k>0. \\
 \frac{ M_{\nu}\left(  k t^{-\nu}\right)  }{t^{\nu}} \div s^{\nu -1}\exp(-k s^\nu),\ 0<\nu<1,\ k>0. \\
 \frac{ W\left(-\nu, \mu;  k t^{-\nu}\right)  }{t^{1-\mu}} \div s^{-\mu}\exp(-k s^\nu),\ 0<\nu<1,\ k>0.
 \end{flalign}
 %%%%%%%%%%%%%%%

 %%%%%%%%%%%%%
 %  Section
 %%%%%%%%%%%%%
 \section{Proofs}\label{sec:proofs}

 %%%%%%%%%%%%%
 %  Section
 %%%%%%%%%%%%%
 \subsection{Proof of Theorem  \ref{prop:raygreen1}}
 \begin{proof}
 	Observe that
 	\[
 	\frac{\partial^k}{\partial x_1^k} \wrighta{ -\frac{\beta}{2}}{ 1} {- \frac{|x_1|}{ t^{\beta /2} \sqrt{D} } } =  \frac{\left( -1\right)^k }{(2 \sqrt{D})^k  t^{\frac{k \beta}{2}} }  \wrighta{ -\frac{\beta}{2}} {1 - \frac{k \beta}{2}} {- \frac{|x_1|}{ t^{\beta /2} \sqrt{D} } }
 	\]
 	Since $n$ is restricted to the unit sphere suppose that $n=e_1$. 
 	Then 
 	\[
 	G_\mathbf{n}^{k}(t) = \frac{(-1)^k}{2 ^k  } \frac{\partial^k}{\partial x_1^k}
 	 \wrighta{ -\frac{\beta}{2}} {1 - \frac{k \beta}{2}} {- \frac{|x_1|}{ t^{\beta /2} \sqrt{D} } }
 	 = \frac{(-1)^k}{2 ^k  } \left( e_1 \cdot \nabla \right) ^k   \wrighta{ -\frac{\beta}{2}}{ 1} {- \frac{|x_1|}{ t^{\beta /2} \sqrt{D} } }.
 	\]
 \end{proof}
 %%%%%%%%%%
 
 %%%%%%%%%%%%%
%  Section
%%%%%%%%%%%%%
\subsection{Proof of Theorem \ref{th:imwright}}

 \begin{proof}
 	The proof technique follows \cite{Luchko2008}. 
	The Wright function is represented by the Hankel integral
	\[
	\wrighta{ a}{ b}{ - z}	= \frac{1}{2 \pi i}\int_{Ha} Ker(\xi) d\xi
	\]
	with kernel
	\[
	Ker(\xi)=\frac{{{ e}^{\xi-\frac{z}{{{\xi}^{a}}}}}}{{{\xi}^{b}}}, a>-1, z>0
	\]
	The contour is depicted in Fig. \ref{fig:hankel}.
	The integral can be split in three parts
	\[
	\int_{Ha} Ker(\xi) d\xi = \int_{AB} Ker(\xi) d \xi + \int_{BCD} Ker(\xi) d \xi + \int_{DE} Ker(\xi) d \xi
	\]
	Therefore, the residue for $b=1$ is given by the limit
	\[
	\llim{\xi}{0}{\xi^{1-b}  { e}^{\xi-\frac{z}{{{\xi}^{a}}}} } = \llim{\xi}{0}{}  { e}^{ -\frac{z}{{{\xi}^{a}}}} 
	\]
	Therefore, for $a<0$ $Res [Ker(\xi)]= \llim{\xi}{0}{  { e}^{ -\frac{z}{{{\xi}^{a}}}}} =0$
	while for $a>0$ $Res [Ker(\xi)]=\llim{\xi}{0}{  { e}^{ -\frac{z}{{{\xi}^{a}}}}} =1$.
	Therefore, in both cases the residue can be neglected by a suitable normalization. 
	Accordingly we can put $ \int_{BCD} Ker(\xi) d \xi = 0 $.
	Along the ray AB $\xi=r e^{i \delta}$ the kernel becomes
	\[
	Ker_A =   \frac{{{e}^{{{  e}^{ i \delta}} r-\frac{z}{{{\left( {{  e}^{  i \delta}} r\right) }^{a}}}}}}{{{\left( {{  e}^{  i \delta}} r\right) }^{b}}}
	\]
	Along the ray DE $\xi=r e^{ - i \delta}$ the kernel becomes
	\[
	Ker_B =   \frac{{{e}^{{{  e}^{- i \delta}} r-\frac{z}{{{\left( {{  e}^{ - i \delta}} r\right) }^{a}}}}}}{{{\left( {{  e}^{  -i \delta}} r\right) }^{b}}}
	\]
	Therefore,
	\[
	Ker_A-Ker_B = \frac{2\,  i \, {{  e}^{\cos{(\delta)} r-\frac{\cos{\left( a \delta\right) } z}{{{r}^{a}}}}}\, \sin{\left( \frac{\sin{\left( a \delta\right) } z}{{{r}^{a}}}+\sin{(\delta)} r-b \delta\right) }}{{{r}^{b}}}
	\]
	Therefore,
	\[
	\llim{\delta}{\pi}{}\frac{1}{2 \pi i} \int_{0}^{\infty} (	Ker_A-Ker_B ) dr= 
	\frac{1}{ \pi }\int_{0}^{\infty} \frac{ {{ e}^{-\frac{\cos{\left(  {\pi}  a\right) } z}{{{r}^{a}}}-r}} }{{{r}^{b}}} \, \sin{\left( \frac{\sin{\left(  {\pi}  a\right) } z}{{{r}^{a}}}- {\pi}  b\right) } dr
	\]
	So that for $b=1$
	\[
	K(a, r) := \frac{{{  e}^{\frac{-\cos{\left( \ensuremath{\pi}  a\right) } z}{{{r}^{a}}}-r}}\, \sin{\left( \frac{\sin{\left( \ensuremath{\pi}  a\right) } z}{{{r}^{a}}}\right) }}{r}
	\]
	Further, for z=0 $K(a,r)=0$.
	
	For $a=-1/2$
	\[
	K\left( -\frac{1}{2}, r\right)  = -\frac{{{  e}^{-r}}\, \sin{\left( \sqrt{r} z\right) }}{r}
	\]
	which, after change of variables the integral becomes\footnote{ \url{http://functions.wolfram.com/GammaBetaErf/Erf/07/01/01/}}
	\[
	-\frac{1}{\pi} \int_{0}^{\infty }{\left. \frac{{{  e}^{-{{y}^{2}}}}\, \sin{\left( y z\right) }}{y}dy\right.} = - \frac{\erf{z/2}}{2}
	\]
	Therefore, 
	\[
	W_I \left( -\frac{1}{2}, -z\right) =\frac{1}{2} - \frac{1}{\pi} \int_{0}^{\infty} 	K\left( -\frac{1}{2}, r\right)  dr,
	\]
	and by the continuous dependence on the b parameter the result follows.
	The case for $z<0$ follows from the symmetry of the Green's function.
\end{proof}
%%%%%%%%%%%%%%%%%%

%%%%%%%%%%%%%%%%
%  REFERENCES
%%%%%%%%%%%%%%%%%
\bibliographystyle{plain}  % Style BST file
\bibliography{biodiffusion1}
\end{document}